\documentclass{article}
\usepackage{arxiv}
\usepackage{amsfonts,amssymb,amsmath,amsthm,latexsym}
\usepackage{graphicx}
\usepackage{algorithm,algorithmic}
\usepackage{multirow}
\usepackage{epstopdf}
\usepackage[labelformat=simple]{subfig}

\usepackage{listings}%
\usepackage{xcolor}
\usepackage{bm}        
\graphicspath{ {./fig/} }

\def\bu{\mathbf{u}}

\def\bg{\mathbf{g}}
\def\bR{\mathbb{R}}

\usepackage{hyperref}

\makeindex

\title{Predicting periodic solutions of cyclic Lotka-Volterra equations with dynamic mode decomposition}

\stitle{Predicting periodic solutions of cyclic LV equations with DMD}

\author{
	\href{https://orcid.org/0000-0003-1037-5431}{\includegraphics[scale=0.06]{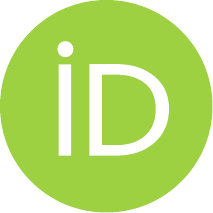}\hspace{1mm}B\"ulent Karas\"ozen} \\
     Department of Mathematics\\
     Middle East Technical University, Ankara-T\"urkiye\\
     \texttt{bulent@metu.edu.tr}\\
       \And
  \href{https://orcid.org/0000-0001-5262-063X}{\includegraphics[scale=0.06]{orcid.pdf}\hspace{1mm}Murat Uzunca} \\
   Department of Mathematics\\
	Sinop University, Sinop-T\"urkiye\\
     \texttt{muzunca@sinop.edu.tr}
}

\begin{document}

\maketitle

\begin{abstract}
Cyclic Lotka-Volterra systems in form of ODEs and PDEs are solved by linearly implicit Kahan's method, that preserves the quadratic Poisson bracket, phase space volume, Hamiltonian ad Casimirs. The solutions are predicted applying dynamic mode decomposition (DMD). Numerical results show that the extended DMD with quadratic dictionaries and Hankel  DMD with the delay embeddings,  predict the periodic solutions with high accuracy, whereas the standard DMD fails. The Hamiltonians and Casimirs are also preserved accurately  by the extend DMD and Hankel  DMD.
\end{abstract}


\section{Introduction}

The Lotka-Volterra (LV) equations mimic the dynamics of interacting species such as predator-prey systems, predicting density oscillations. 
For spatially inhomogeneous situations, LV equations are straightforwardly generalized to diffusion-reaction equations, widely applied to more complex ecological processes.  
The LV equation in Poisson form is solved by the average vector field (AVF) method \cite{Cohen11} and  by line integral method \cite{Brugnano12}. 
It was shown that the Hamiltonian and the Casimir are preserved in long-term. 
The Volterra lattice is solved by Lobatto IIA-B method \cite{Karasozen06} preserving  the Poisson structure, the Hamiltonian and the Casimirs.

Dynamic mode decomposition (DMD)  is simple to implement and yields a linear representation of the nonlinear dynamics, which is well-suited for prediction, analysis, and control. 
DMD is also related to Koopman operator theory \cite{Brunton22}, which offers a linear perspective on inherently nonlinear systems by operating in an infinite-dimensional space of observable functions. 
Several variants of the DMD algorithm have been developed, such as the standard DMD, exact DMD etc. \cite{Tu14}. 
In the extended DMD (EDMD) \cite{Williams15}, dictionaries of observable
functions are used to learn the finite-dimensional projection of the Koopman operator.  
Hankel with delay embeddings (HDMD) \cite{Champion19,Kamb20} yields more accurate predictions than  DMD, which is particularly effective for ergodic systems that exhibit low-dimensional attractors. 
Prediction accuracy of the DMD was investigated in \cite{Lu20}. 
Time series of  linear and nonlinear systems are predicted by EDMD \cite{Zinihi26}, HDMD \cite{Tirunagari17} and by regularized DMD \cite{Xie24}.

In this paper we solve the cyclic LV equations with periodic solutions in form ordinary differential equation (ODE) and partial differential equation (PDE)  with the second order linearly implicit Kahan's method  \cite{Celledoni13,Kahan93,Kahan97}.  
Kahan's method is designed for solving a general linear-quadratic system.  It preserves the quadratic Poisson structure, the linear Hamiltonian,  quadratic and cubic Casimirs and the phase volume. 
We perform predictions with the DMD, EDMD, HDMD and extended delayed DMD (EHDMD). 
Numerical results for three and four component cyclic LV equations show that the prediction accuracy increases by EDMD and by the number of shifting parameter of HDMD. 
In the case of the diffusive LV equation, the same prediction behavior is observed; the solutions are confined on a closed curve and the total mass  is preserved for the EDMD, HDMD and EHDMD, whereas the  DMD fails to predict periodic solutions.
presented.

The outline of the paper is as follows. 
In Sections~\ref{sec:LV}-\ref{sec:difLV}, we introduce the cyclic LV systems and diffusive LV equation, respectively. 
In Section~\ref{sec:kahan}, Kahan's method and its discrete integrability properties are summarized. 
In Section~\ref{sec:dmd}, the dynamic mode decomposition, Koopman operator and time-delayed embeddings are introduced. Numerical results for three and four dimensional LV systems and three component diffusive LV equation are presented in Section~\ref{sec:num}.   
The paper end with some conclusions in Section~\ref{sec:conc}.

\section{Hamiltonian Lotka-Volterra equations}
\label{sec:LV}

The $n$-dimensional LV system is described by the following system of differential equations:
\begin{equation}\label{lvs}
  \dot u_i = b_i u_i + \sum_{j=1}^{n} a_{i,j} u_i u_j, \ \ i=1,2, \dots , n \; .
\end{equation}
where $u_i=u_i(t)$ denotes the concentration or population of the $i$-th species, $b_i$ is the rate constant of the $i$-th species, and $a_{ij}$ are the interaction coefficients between the $i$-th $j$-th species.
In this paper, we consider  LV system in skew-symmetric form , i.e., $a_{ij}=-a_{ji}$ for $1\leqslant i,j\leqslant n$ and $b_i=0$ \cite{Vanhaecke22}.
The best known example of an integrable LV equation is the $n$-particle periodic or cyclic Kac-van Moerbeke (KM) system \cite{Kac75}, given by
\begin{align}\label{KM}
  \dot{u}_{i}&=u_{i}(u_{i+1}-u_{i-1}), & i=1,\ldots,n.
\end{align}

Defining the solution vector $\bu (t) = (u_1(t), u_2(t), \cdots,u_n(t))^T: [0,T] \mapsto \bR^{n}$, the KM system \eqref{KM} can be written as
\begin{equation}\label{cKM}
\dot{\bu} = \bu\odot (A\bu),
\end{equation}
where $\odot$ denotes the element-wise multiplication (Hadamard), and $A= (a_{ij})$ is the interaction matrix defined by
\begin{equation}\label{KMA}
 A=\left(
  \begin{array}{ccccrc}
    0&1& 0 &\cdots& 0 &-1\\
     -1 &0&1&0      &\dots   & 0\\
     0 & -1 &0   &      &   &\vdots\\
    \vdots&&\ddots&\ddots& &0\\
       0 &  & & & 0&1\\
    1&0&\cdots &\cdots&-1&0
  \end{array}
  \right)\,,
\end{equation}
with the indices are periodic modulo $n$, i.e., $u_{n+1}=u_1$ and $u_0=u_n$.  
The KM system  \eqref{KM} also possess a Hamiltonian structure 
$$
\dot{ \bu} = J(\bu) H(\bu),
$$
with the quadratic Poisson brackets $\{u_i,u_j\}:=u_iu_j$, $1\leq i,j\leq n$. Here, $J(\bu)$
is  the  skew-symmetric matrix, given by
$$
 J(\bu)=\left(
  \begin{array}{ccccr}
    0&u_1u_2& \cdots &\cdots& -u_1u_n\\
         \vdots&\ddots&\ddots& &0\\
     \cdots & -u_{i-1}u_i &0   & u_iu_{i+1}     &   \cdots\\
    \vdots&\vdots&\vdots&\ddots &\vdots\\
    u_1u_n&\cdots &\cdots&-u_{n-1}u_n&0
  \end{array}
  \right)\,,
$$
and $H(\bu) = \sum_{i=1}^n u_i$  the linear Hamiltonian. 
It is also known as  the Volterra lattice, for the integrable discretization of the Korteweg de-Vries (KdV) equation \cite{Kac75,Suris99} and of the inviscid Burger's equation \cite{Kupershmidt97}. 
The Volterra lattice is also bi-Hamiltonian with the cubic Poisson bracket \cite{Karasozen06}. 
The $n$-dimensional extension of \eqref{KM} as integrable discretization of the Korteweg de-Vries equation was integrated with a Poisson structure preserving integrator in \cite{Karasozen06}.
The rank of the interaction matrix $A$ in \eqref{KMA} is $n-1$ when $n$ is odd and $n-2$ otherwise. 
In the first case, when $n$ is odd, the product $u_1u_2\dots x_n$ is a Casimir function. In the other case, both the products $u_1u_3\dots u_{n-1}$ and $u_2u_4\dots u_n$ are Casimir functions. 
An additional $(n-1)/2$ independent polynomial first integrals (including the Hamiltonian), in involution, are constructed from a Lax equation.  
This accounts for the Liouville integrability of KM system  \eqref{KM} \cite{Vanhaecke22}.

\section{Diffusive Lotka-Volterra equations}
\label{sec:difLV}

In this section, we consider the $n$-dimensional LV system including diffusion terms, on the time period $[0,T]$ and on a smooth spatial domain $\Omega\subset\bR^2$ \cite{Suzuki15} 
\begin{equation} \label{lvsd}
\frac{\partial u_i}{\partial t}  = d_i \Delta u_i  + u_i\left ( \sum_{j=1}^n a_{ij}u_j\right), \quad i=1,\ldots,n,
\end{equation}
with the no-flux Neumann boundary conditions  
and with the initial conditions
$$
\qquad u_i(x,y,0) = U_{i}(x,y) \ge 0,
$$
where $d_i > 0 $ are the diffusion rates, $a_{ij}$ are the inter-specific interaction coefficients between the species as in \eqref{lvs}, and $\Delta$ is the second order Laplace operator.
The dynamics of the system \eqref{lvsd} is controlled by the ODE part and the solution becomes spatially homogeneous in form of patterns oscillating in time. 
We consider the following three species diffusive LV system \cite{Suzuki15}
\begin{equation}\label{lvs2}
\begin{aligned}
\frac{\partial u_1}{\partial t} & =  d_1 \Delta u_1 + u_1( u_2- u_3), \\
\frac{\partial u_2}{\partial t} & =  d_2 \Delta u_2  +  u_2(u_3 - u_1), \\
\frac{\partial  u_3}{\partial t}& =  d_3 \Delta u_3  +  u_3( u_1 - u_2). 
\end{aligned}
\end{equation}

The total mass $M(u)$ of the system \eqref{lvs2}  is conserved as
\begin{align} \label{tm}
\frac{d}{dt} M(u) &=0 , & M(u)=\sum_{j=1}^3 \int_ \Omega u_jd\Omega.
\end{align} 

For two component diffusive LV systems, there exits a logarithmic Hamiltonian \cite{Suzuki11, Settani16}, and due to non-invertibility of the coefficient matrix  $A$ in \eqref{KMA},  the total mass $M(u)$ is the only conserved quantity.

Finite-difference space discretization of the diffusive LV system \eqref{lvs2} leads to the  coupled linear-quadratic ODE system of the following form
\begin{equation}\label{dlvs}
\begin{aligned}
\frac{d \bu_1}{d t} &= d_1D\bu_1 + \bu_1\odot(\bu_2 - \bu_3), \\
\frac{d \bu_2}{d t} &= d_2D\bu_2 + \bu_2\odot(\bu_3 - \bu_1), \\
\frac{d \bu_3}{d t} &= d_3D\bu_3 + \bu_3\odot(\bu_1 - \bu_2),
\end{aligned}
\end{equation}
where $\bu_i(t): [0,T] \mapsto \bR^{m^2}$ are the semi-discrete approximations to the exact solutions $u_i(x,y,t)$, given by the ordering
$$
\bu_i(t) = (u_{i}(x_1,y_1,t),\ldots,u_{i}(x_m,y_1,t),\ldots,u_{i}(x_m,y_1,t),\ldots,  u_{i}(x_m,y_{m},t))^T,
$$
at the uniform spatial grid nodes $(x_i,y_j)\in\Omega$, $1\leq i,j \leq m$.
The matrix $D\in\bR^{m^2\times m^2}$ in the system \eqref{dlvs} stands for the discrete Laplace operator, which can be computed from the one-dimensional spatial setting utilizing the Kronecker product $\otimes$. 
Let the matrix $D_{(1)}\in \mathbb{R}^{m\times m}$ be the matrix of the one-dimensional second-order differential operator $\partial^2/ \partial_{xx}$ under homogeneous Neumann boundary conditions with the spatial mesh size $\Delta x$.
Then, the matrix $D$ of the discrete Laplace operator on the given two-dimensional spatial grid is computed by
$$
D = D_{(1)}\otimes I_m + I_m \otimes D_{(1)} \in  \bR^{m^2\times m^2},
$$
where $I_m\in \mathbb{R}^{m\times m}$ is the $m$-dimensional identity matrix. 
For the solution vector of the system, we introduce $\bu(t)=(\bu_1(t)^T,\bu_2(t)^T,\bu_3(t)^T)^T: [0,T] \mapsto \bR^{3m^2}$. Then, the system \eqref{dlvs} can be written in the compact form
\begin{equation}\label{cdlvs}
\frac{d \bu}{d t} = L\bu + Q(\bu),
\end{equation}
where the constant matrix $L\in\bR^{3m^2\times m^2}$ and the quadratic vector $Q\in\bR^{3m^2}$ are defined by
$$
L = 
\begin{pmatrix}
d_1D &  &  \\
 & d_2D & \\
&& d_3D
\end{pmatrix}, \qquad 
\begin{pmatrix}
\bu_1\odot(\bu_2   -\bu_3)  \\
\bu_2\odot(\bu_3   -\bu_1) \\
\bu_3\odot(\bu_1   -\bu_2)
\end{pmatrix}.
$$

By the given discretization setting, the discrete form of the total mass $M(u)$ in \eqref{tm} is given by
$$
M(\bu (t)) = \Delta x \Delta y \sum _{l=1}^3 \sum _{i,j=1}^m  u_l (x_i, y_j,t).
$$

In \cite{Settani16}, it has been shown  that associated Hamiltonian   of  the two component diffusive  LV system is preserved by the symplectic scheme, and the spatially average solutions are periodic in time.

\section{Kahan's discretization}  
\label{sec:kahan}

The Kahan's method is designed for solving a general linear-quadratic system of the form 
$$
\dot{\bu} = f({\bu}) =  Q({\bu}) + L{\bu},
$$
leading to the scheme
\begin{equation}\label{kahan}
\frac{\bu^{k+1} - \bu^n}{\Delta t} = \widetilde{Q}(\bu^k,\bu^{k+1}) + \frac{1}{2}L(\bu^k + \bu^{k+1}),
\end{equation}
where, $\Delta t$ is the time-step size and the symmetric bilinear form $\widetilde{Q}(\cdot ,\cdot )$ is obtained by the polarization of the quadratic vector field $Q(\cdot)$ 
\begin{equation}\label{polar}
\widetilde{Q}(\bu^k,\bu^{k+1}) := \frac{1}{2}\left( Q(\bu^k+\bu^{k+1}) - Q(\bu^k) -  Q(\bu^{k+1}) \right).
\end{equation}

The solution  $\bu^{k+1}$  of  \eqref{kahan} can be computed by solving a single linear system of equations \cite{Celledoni13}
$$
\frac{\bu^{k+1} -\bu^k}{\Delta t} =\left ( I -\frac{\Delta t}{2} f'(\bu^k)\right )^{-1} f(\bu^k),
$$
where $f'$ denotes the Jacobian matrix of $f$. 
For the approximation $u_i^k:=u_i(t_k)$ at time $t_k$, Kahan's method for the KM system \eqref{KM} yields for $i=1,\ldots ,n$
\begin{align*}
\frac{u_i^{k+1} - u_i^k}{\Delta t } &= \frac{1}{2}\sum_{i=1}^{n}  \left (u_i^{k+1}u_{i+1}^{k} + u_i^{k+1} u_{i+1}^{k}\right ) -  \left (u_i^{k+1}u_{i-1}^{k} + u_i^{k+1} u_{i-1}^{k}\right )  , 
  &k=1,\ldots, 
\end{align*}
which preserves the quadratic Poisson bracket 
$
\{u_i^{k},u_j^{k}\}=u_i^{k}u_j^{k}
$, i.e., it is a Poisson integrator \cite{Kouloukas16,Vanhaecke22}. 
For the approximation $\bu^k:=\bu(t_k)$ at time $t_k$, Kahan's method for the KM system \eqref{cKM} in compact form, yields
\begin{align*}
\frac{\bu^{k+1} - \bu^k}{\Delta t } &= \frac{1}{2} \left ( \bu^{k+1}\odot (A\bu^{k}) + \bu^{k}\odot (A\bu^{k+1}) \right) .
\end{align*}

On the other hand, Kahan's method for the diffusive LV system \eqref{cdlvs} yields
\begin{align*}
\frac{\bu^{k+1} - \bu^k}{\Delta t } &= \frac{1}{2} L\left ( \bu^{k+1} + \bu^{k} \right) + \widetilde{Q}(\bu^k,\bu^{k+1}),
\end{align*}
where $\widetilde{Q}$ is the symmetric bilinear form  defined in \eqref{polar}.

A dynamical system $\dot{\bu} =f(\bu)$ is called reversible if $f(\rho\bu)=-\rho f (\bu)$,  where $\rho$ is an invertible linear transformation in the phase space. 
The KM system \eqref{KM} is reversible with $\rho = -1$.  
The time integrator  $\bu^{k+1} =\Phi (\bu^k)$ is symmetric or time-reversible if it holds $\Phi_{\Delta t} =\Phi^{-1}_{-\Delta t}$  \cite{Hairer06}. 
Kahan's method is time-reversible \cite{Celledoni12}, 
$$
\frac{\bu^{k+1} -\bu^k}{\Delta t} = \left ( I +\frac{\Delta t}{2} f'(\bu^{k+1})\right )^{-1} f(\bu^{k+1}),
$$
and hence it is a second-order integrator \cite{Celledoni13}.

The vector field of the conservative part of the KM equation \eqref{KM} is divergence-free: $\nabla f(\bu) = 0$, hence, the phase-space volume is preserved.
For Hamiltonian systems with quadratic vector fields, Kahan's method preserves the modified
volume and modified Hamiltonians \cite{Celledoni14,Iserles15}.

\section{Dynamic mode decomposition}
\label{sec:dmd}

In the sequel, we consider the data vectors $\{\bu^1,\ldots ,\bu^{n_t}\}$ as the solution vectors of the KM system \eqref{cKM} or the diffusive LV system \eqref{cdlvs}, where $n_t$ is the number of discrete time instances. Without lost of generality, we set $n_x$ as the length of each data vector $\bu^k$, for easy notation, i.e., $\bu^k\in\bR^{n_x}$ for $k=1,\ldots ,n_t$.

For a constant matrix $K \in \mathbb{R}^{ n_x \times n_x}$, the DMD algorithm estimates a linear relationship 
$$
	U_2  \approx K U_1, 
$$
between the data matrices
$$
    U_1 = 
		\begin{pmatrix}
    \vert&\vert&&\vert\\
   \mathbf{u}^1   & \mathbf{u}^2 & \cdots & \mathbf{u}^{n_t-1}\\
    \vert&\vert&&\vert
    \end{pmatrix}\in \mathbb{R}^{ n_x \times n_t-1}, \qquad 
     U_2 = 
		\begin{pmatrix}
    \vert&\vert&&\vert\\
     \mathbf{u}^2  & \mathbf{u}^3 & \cdots & \mathbf{u}^{n_t}\\
    \vert&\vert&&\vert
    \end{pmatrix}\in \mathbb{R}^{ n_x \times n_t-1}.
$$
The optimal $K$ is found by solving the optimization problem
$$
K = \arg\min\limits_{\hat{K}} \Vert \hat{K}U_1 - U_2 \Vert_F,
$$
where $|| \cdot||_F$ denotes the Frobenius norm defined as
    $\Vert U \Vert_F = \sqrt{\sum_{i }\sum_{j } U_{ij}^2}$.
The least-squares solution to this optimization problem is known to be
\begin{equation} \label{dmds}
   A = U_2 U_1^{\dagger},
\end{equation}
where $U_1^{\dagger}$ is the Moore-Penrose inverse of the matrix $U_1$.  
There are different versions of the DMD algorithm such as standard DMD and exact DMD \cite{Tu14}. 
In this paper, we use standard DMD, given in Algorithm~\ref{dmd}.

\begin{algorithm}
\caption{Standard DMD}\label{dmd}
\begin{algorithmic}[1]
\REQUIRE Solution vectors $\{\bu^1,\ldots,\bu^{n_t}\}$
\ENSURE DMD modes $\{\hat{\phi}_1,\ldots,\hat{\phi}_r\}$ 
\STATE Arrange the data $\{\bu^1,\ldots,\bu^{n_t}\}$ into the matrices
        $$
            X = \begin{bmatrix} \bu^1 & \cdots & \bu^{n_t-1} \end{bmatrix},
            \qquad
            Y = \begin{bmatrix} \bu^2 & \cdots & \bu^{n_t} \end{bmatrix}
        $$
\STATE	Compute the (reduced) SVD of $X$, 
        $X = U \Sigma V^* $, with rank $r$
\STATE Define the matrix $ \tilde A = U^* Y V \Sigma^{-1}\in\bR^{r\times r}$
\STATE Compute eigenvalues and eigenvectors of $\tilde A$, 
        $\tilde A w_j = \lambda_j w_j$
        
\STATE  Compute the DMD mode $\hat{\phi}_j\in\bR^{n_x}$ corresponding to the DMD eigenvalue $\lambda_j$, 
        $ \hat{\phi} = U w_j$     
\end{algorithmic}
\end{algorithm}

Each column of $\Phi:=[\hat{\phi}_1 \; \cdots \; \hat{\phi}_r]\in\bR^{n_x\times r}$ in Algorithm~\ref{dmd} is a DMD mode corresponding to a particular eigenvalue in $\Sigma$. 
With the approximated eigenvalues and eigenvectors of $A$, a solution at some $k$-th time step ($k>n_t$) can be constructed explicitly as
\begin{align}\label{eq:2-7}
\bu_\text{DMD}^{k} &= \Phi\Sigma^{n+1}{\mathbf b},  & k>n_t,
\end{align}
where ${\mathbf b} =\Phi^{-1}\bu^{1}\in\bR^{ r}$ is the vector representing the initial amplitude of each mode. 
The solution at any future time is approximated directly with \eqref{eq:2-7} using only information encapsulated in the first $n_t$ temporal snapshots.

\subsection{Extended dynamic mode decomposition}

The Koopman operator framework provides a powerful alternative to traditional nonlinear analysis by shifting the dynamics from the original state space to a higher-dimensional (possibly infinite-dimensional) space of observables, where, the evolution becomes linear. 
This property allows to use linear operator theory to predict and  control on nonlinear systems.
For a nonlinear dynamic system $\dot{\bu} = f(\bu)$, the Koopman operator $\mathcal K$ is an infinite-dimensional linear operator that acts on all observable functions $g: \mathcal M\to \mathbb C$ so that \cite{Brunton22}
\begin{equation}\label{kop}
\mathcal K g(\bu) = g(f(\bu)).
\end{equation}
For a discrete dynamic system, the discrete-time Koopman operator $\mathcal K_t$ is 
$$
\mathcal K_t g(\bu^{k}) = g(\mathcal K_t(\bu^k)) = g(\bu^{k+1}).
$$
The Koopman operator transforms the finite-dimensional nonlinear problem in the state space into the infinite-dimensional linear problem in the observable space. Since $\mathcal K_t$ is an infinite-dimensional linear operator, it is equipped with infinite eigenvalues $\{\lambda_i\}_{i=1}^{\infty}$ and eigenfunctions $\{\phi_i\}_{i=1}^\infty$. 
In practice, the eigenvalues and eigenfunctions has to be approximated in a finite-dimensional space.

Let $\mathbf g$ denotes a $p \times 1$ vector of observables,
$$
{\mathbf g}(\bu^{k}) = \begin{bmatrix}
g_1(\bu^k)\\
\vdots\\
g_p(\bu^k)
\end{bmatrix},
$$
where $g_j: \mathcal M \to \mathbb C$ is an observable function, with $j =1,\cdots, p$. If the chosen observable $\mathbf g$ is restricted to an invariant subspace spanned by eigenfunctions of the Koopman operator $\mathcal K_t$, then it induces a linear operator $K$ that is finite-dimensional and advances these eigen-observable functions on this subspace \cite{Brunton22}.

Based on \eqref{kop}, the DMD algorithm can be deployed to approximate the eigenvalues and eigenfunctions of $K$ using the collected temporal snapshots in the observable space, which is known as the extended DMD \cite{Williams15}. One can compute the DMD on the  lifted data  matrices:
\begin{align*}
    \bg(U_1) &= \begin{pmatrix}
    \vert&\vert&&\vert\\
   \bg(\mathbf{u}^1)   & \bg(\mathbf{u}^2) & \cdots & \bg(\mathbf{u}^{n_t-1})\\
    \vert&\vert&&\vert
    \end{pmatrix}
    & \bg(U_2) = \begin{pmatrix}
    \vert&\vert&&\vert\\
     \bg(\mathbf{u}^2)  & \bg(\mathbf{u}^3) & \cdots & \bg(\mathbf{u}^{n_t})\\
    \vert&\vert&&\vert
    \end{pmatrix}.
\end{align*}

The future state $\bu_\text{DMD}^{k+1}$ is predicted  as 
$$
\bu_\text{DMD}^{k+1} = \Phi\Sigma^{k+1} {\mathbf b}, \qquad  {\mathbf b}= \Phi^{-1}\bu^1, \qquad k>n_t.
$$
Transform from observables space back to the state space,
$$
\bu_\text{DMD}^k ={\mathbf  g}^{-1}(\bu_\text{DMD}^k).
$$
Judicious selection of the observables is critical to success of the Koopman method.
The form of the KM system \eqref{KM} suggests a set of quadratic observables $g_j=u_ju_{j+1}$, $j=1,\ldots,n$.

\subsection{Hankel dynamic mode decomposition}

Hankel-DMD, introduced by \cite{Arbabi17},  represents a specialized instance of EDMD where the dictionary is constructed through time-delay embedding. 
Given a single trajectory of the observable, $\{g(\bu^1),g(\bu^2),\ldots,g(\bu^{n_t})\}$, the matrices $H_1$ and $H_2$  are given explicitly by the Hankel matrices
$$
H_1=\begin{pmatrix} 
g(\bu^1) & g(\bu^2) &  \cdots & g(\bu^{n_t-q})\\
g(\bu^2) & g(\bu^3) &  \cdots & g(\bu^{n_t-q+1})\\
\vdots & \vdots & \vdots & \vdots  \\
g(\bu^{q}) & g(\bu^{q+1}) &  \cdots & g(\bu^{n_t-1})
\end{pmatrix},\quad
H_2=\begin{pmatrix} 
g(\bu^2) & g(\bu^3) &  \cdots & g(\bu^{n_t-q+1})\\
g(\bu^3) & g(\bu^4) &  \cdots & g(\bu^{n_t-q+2})\\
\vdots & \vdots & \vdots & \vdots  \\
g(\bu^{q+1}) & g(\bu^{q+2}) &  \cdots & g(\bu^{n_t})
\end{pmatrix}.
$$

This approach is particularly effective for ergodic systems that exhibit low-dimensional attractors. 
The Hankel alternative view of Koopman (HAVOK) method \cite{Champion19,Kamb20},  based on
time-delay embedding coordinates, can be used to obtain a linear model that nearly perfectly captures the dynamics of nonlinear quasi-periodic systems on the attractor.

Without time-shifting the data, the DMD approximation does not capture the correct
complex eigenvalue pairs associated with the periodic (Fourier) time dynamics. 
Thus $U_1$ and $U_2$  are linearly consistent if and only if the nullspace of $U_2$ contains the nullspace of $U_1$. 
When the data are not linearly consistent, the equation \eqref{dmds} is not exactly satisfied. 
The Koopman analogy can break down and so the DMD analysis \cite{Tu14}.
The time-delay embedding can significantly improve upon the DMD algorithm for producing an approximate dynamical system for forecasting.  
For a short time-delay embedding (small $q$), the time-shifted data can provide a more accurate
assessment of the true rank of the underlying system.
For long time-delay embeddings (large $q$), the nonlinear dynamics can be made to be approximately linear, thus providing an approximation to the Koopman operator and a linear
reduced order modeling \cite{Champion19,Kamb20}.

\section{Numerical results}
\label{sec:num}

In this section, we present numerical results for the KM system and the diffusive LV system.
Both the systems are integrated with the time-step size $\Delta t = 0.01$.  
For a time-dependent Hamiltonian $H(t)$ and a Casimir $C(t)$, the accuracy of the preservation is measured by the relative errors defined by
$$
R_H=\frac{H(t) - H(0)}{H(0)}, \qquad R_C = \frac{C(t) - C_0}{C_0}.
$$
In case of solution accuracy, we use the $L^2$-error between the numerical solution and the DMD approximations. We also give the speed-up factors computed by the ratio of the time needed to obtain the numerical solution over the time needed to compute the DMD approximations.

\subsection{3D KM system} \label{3KM}

A well known cyclic three dimensional LV system is given  as  \cite{Damianou09} 
\begin{align*}
\dot{u}_1 & =  u_1(u_2 - u_3),   \\
\dot{u}_2 & =  u_2(u_3 - u_1),  \\
\dot{u}_3 & =  u_3(u_1 - u_2),  
\end{align*}
represents
Belousov-Zhabotinskii system of three reactants \cite[pp.16]{Perthame15} with the Hamiltonian  $H(t) = u_1(t)+ u_2(t)+ u_3(t)$, and with the Casimir $ C(t)= u_1(t)u_2(t)u_3(t)$. 
For the simulation, we take the initial conditions 
$u_1(0) = 0.5$,  $u_2(0) = 0.5$ and $u_3(0) = 1.5$.
As the observables, $u_1u_2$, $u_1u_3$ and $u_2u_3$ are used in the EDMD.

\begin{figure}[ht]
\centering
\hspace*{-1cm}\includegraphics[scale=0.3]{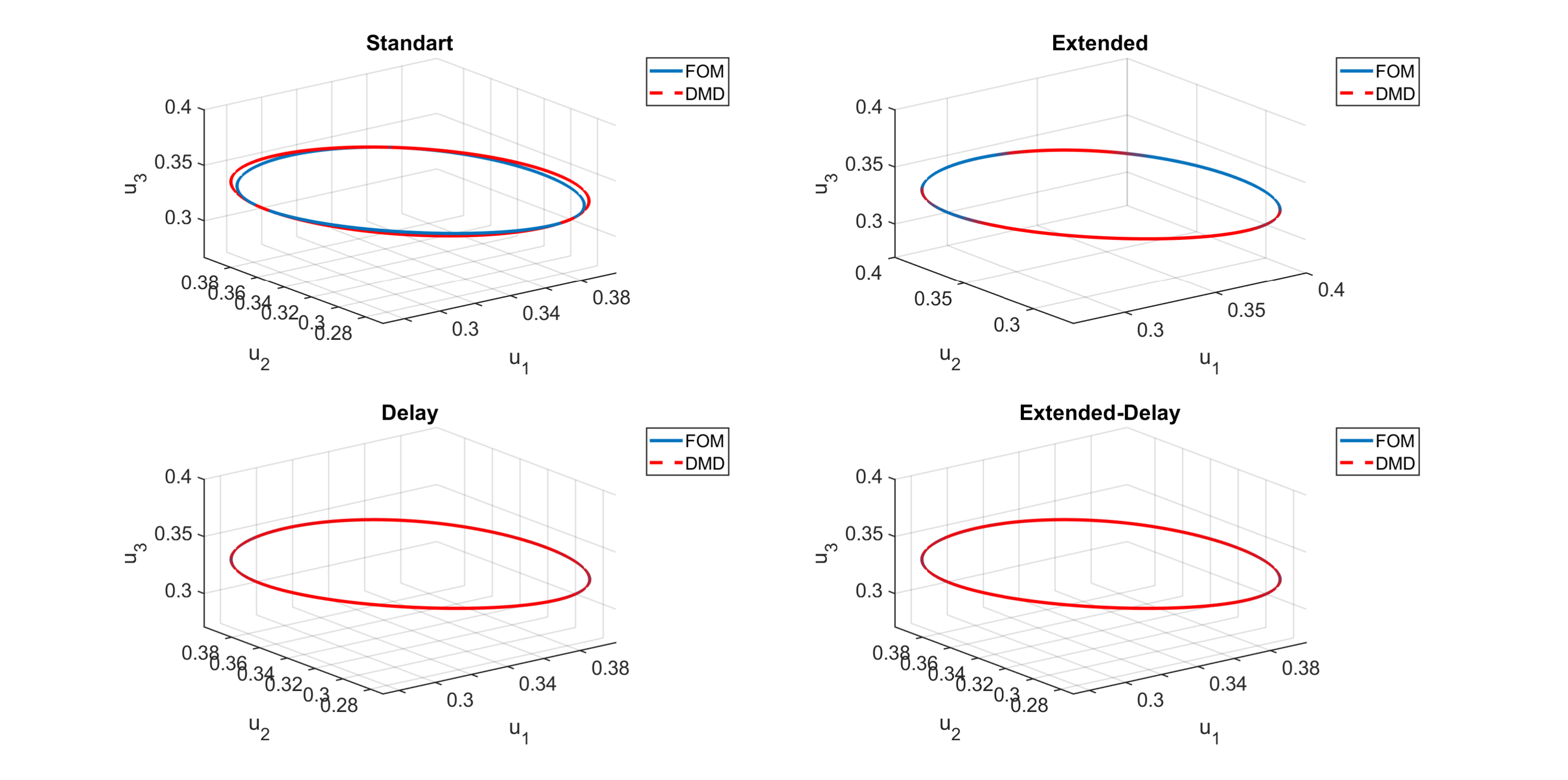}
\caption{Phase portraits of KM system}
\label{lv3phase}
\end{figure}

\begin{figure}[ht]
\centering
\hspace*{-1cm}\includegraphics[scale=0.3]{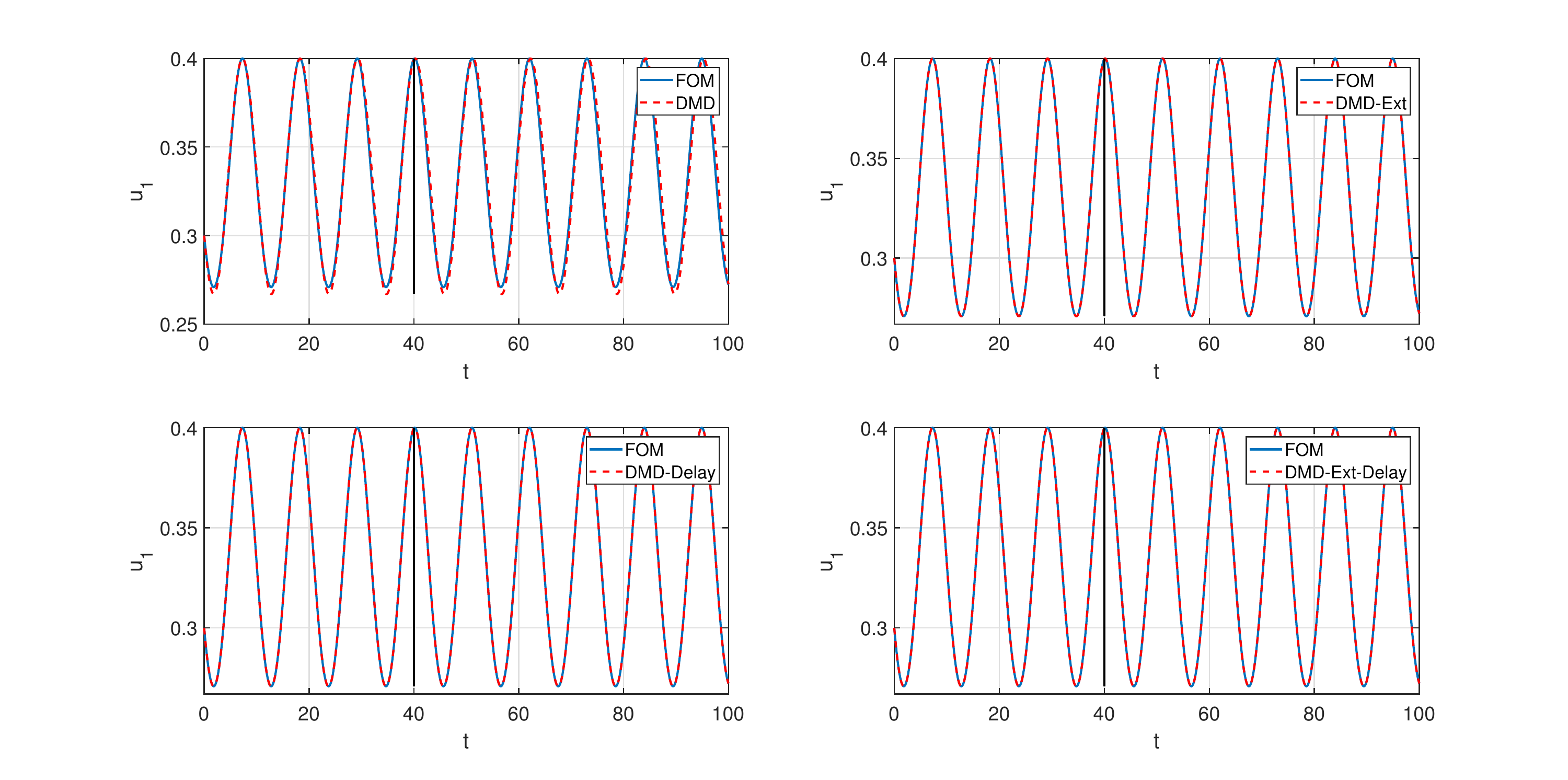}
\caption{Solution profiles of KM system}
\label{lv3sol}
\end{figure}

\begin{figure}[ht]
\centering
\hspace*{-2cm}\includegraphics[scale=0.32]{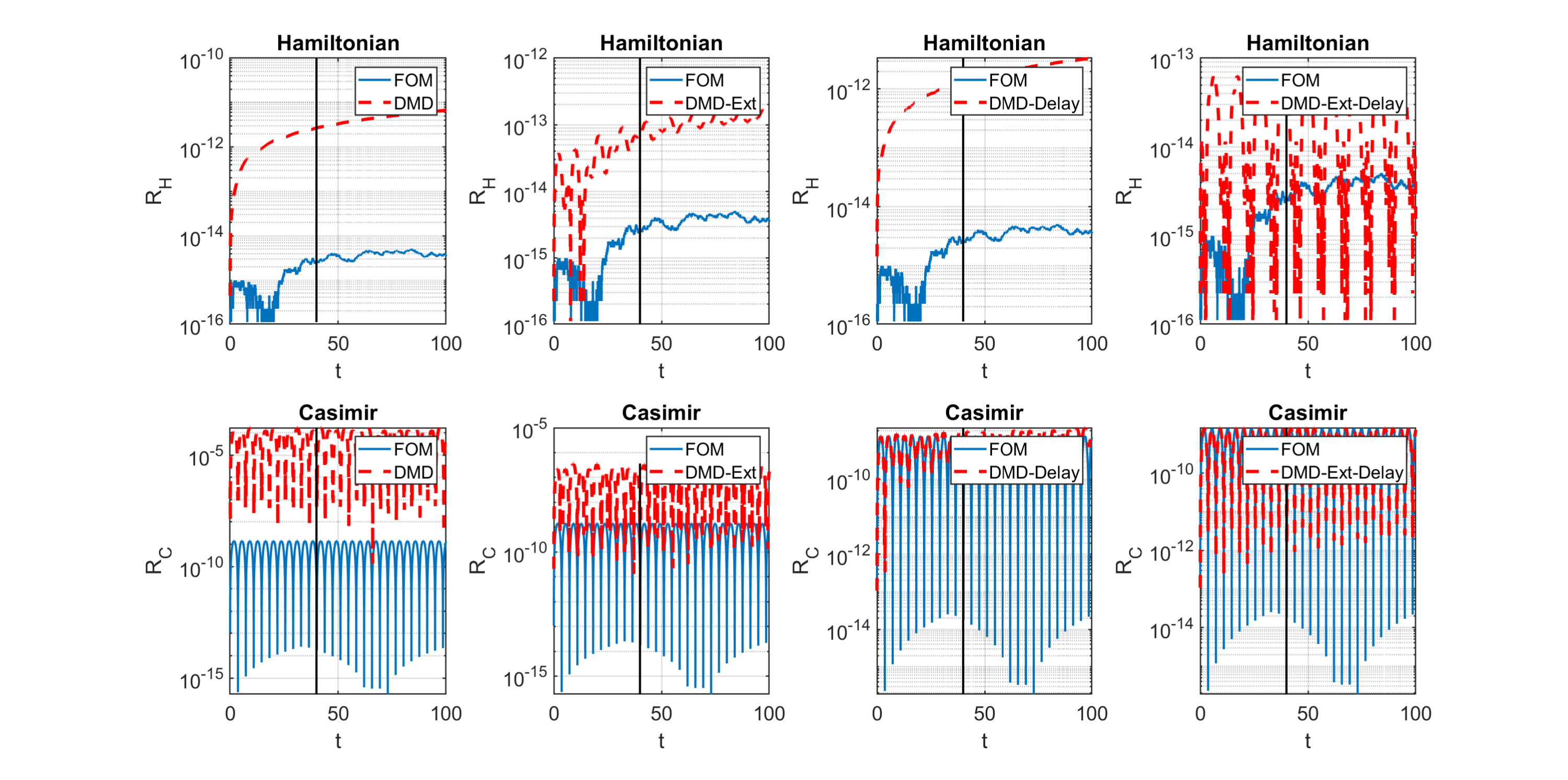}
\caption{Hamiltonian \&  Casimir errors for KM system}
\label{lv3int}
\end{figure}

\begin{table}[ht]
\centering
\caption{Prediction errors for KM system (50-100)}\label{error3}%
\begin{tabular}{@{}lccc@{}}
\hline
Method & r & q & $L_2$ error \\ 
\hline
DMD & 3 & 1 & 9.23e-01 \\
EDMD & 6 & 1 & 7.69e-04 \\
HDMD & 15 & 5 & 2.30e-06\\
HDDMD & 30  & 5  & 9.01e-08\\
DDMD & 30 & 1 & 3.29e-07 \\
EHDMD & 60  & 10  & 5.31e-08 \\
\hline
\end{tabular}
\end{table}

Periodic solutions are accurately predicted with EDMD and HDMD in Figures \ref{lv3phase}-\ref{lv3sol}.  
The linear Hamiltonian is preserved up to machine accuracy and the Casimir errors oscillate without showing any drift over time in Figure \ref{lv3int}.
EHDMD can improve the prediction accuracy as shown in Table~\ref{error3}.

\subsection{4D Volterra lattice} 
\label{4Dvolt}

We consider the four dimensional Volterra lattice
\begin{align*}
\dot{u}_1 & =  u_1(u_2 - u_4),  \\
\dot{u}_2 & =  u_2(u_3 - u_1),  \\
\dot{u}_3 & =  u_3(u_4 - u_2),  \\
\dot{u}_4 & =  u_4(u_1 - u_3),  
\end{align*}
with the Hamiltonian  $H(t) = u_1(t)+ u_2(t)+ u_3(t) +u_4(t)$, and with the Casimirs $ C_1(t)= u_1(t)u_3(t)$  and $C_2(t) = u_2(t)u_4(t)$.
The initial conditions are
$u_1(0) = 0.5$,  $u_2(0) = 0.5$ and $u_3(0) = 1.5$.
As observables, $u_1u_2$, $u_1u_3$, $u_1u_4$, $u_2u_3$, $u_2u_4$ and $u_3u_4 $ are used in the EDMD.

\begin{figure}[ht]
\centering
\hspace*{-1.8cm}\includegraphics[width=1.2\textwidth]{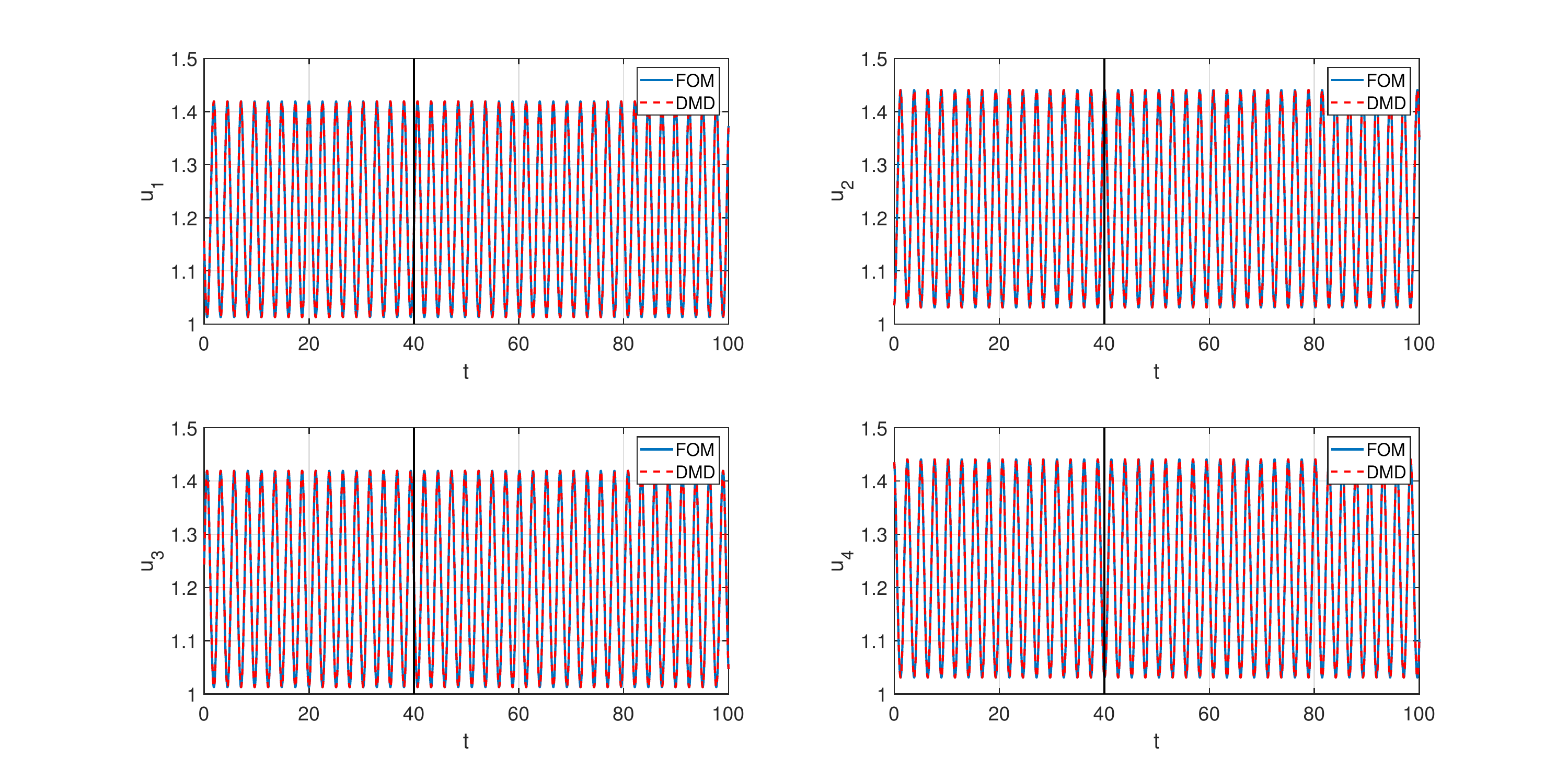}
\caption{Numerical and DEDMD solution profiles of 4D Volterra lattice}\label{lv4sol}
\end{figure}

\begin{figure}[ht]
\centering
\hspace*{-1cm}\includegraphics[width=1.1\textwidth]{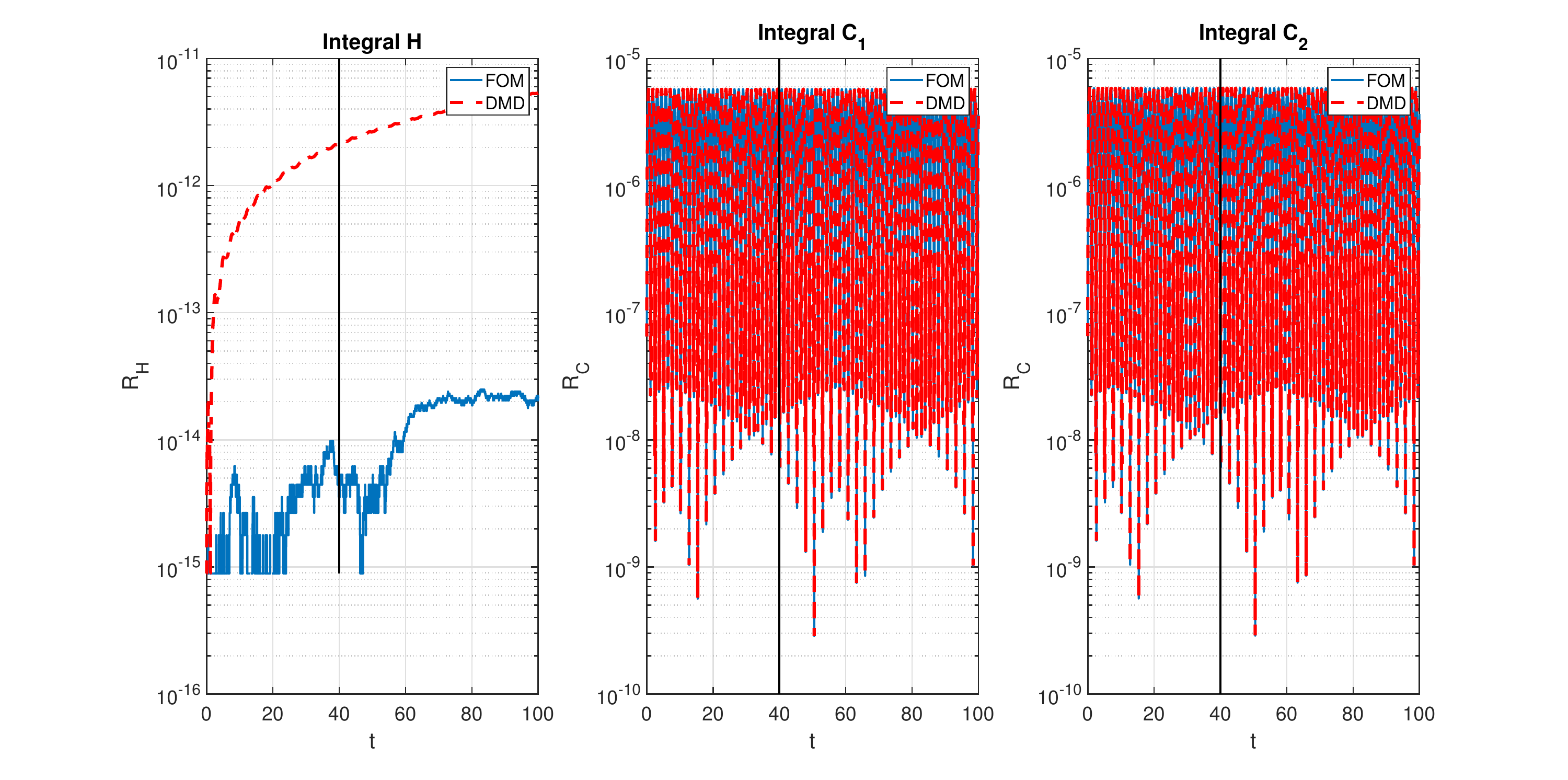}
\caption{Hamiltonian \&  Casimir errors for 4D Volterra lattice by numerical and DEDMD solutions}
\label{lv4ham}
\end{figure}

\begin{table}[ht]
\caption{Prediction errors for 4D Volterra lattice (40-100)}\label{error4}%
\begin{tabular}{@{}lccc@{}}
\hline
Method & r & q & $L^2$-error \\ 
\hline
DMD & 4& 1 & 2.28e+00 \\
EDMD & 10 & 1 & 5.22e-07 \\
HDMD & 40 &10 &1.72e-08\\
HDDMD & 100 & 10 & 2.90e-10  \\
\hline
\end{tabular}
\end{table}

As for the three species LV system, HDMD can improve the prediction accuracy in Figures \ref{lv4sol}-\ref{lv4ham}, and in Table \ref{error4}.

\subsection{Diffusive LV system with periodic solutions}

We consider the diffusive LV system \eqref{lvs2} on the spatial domain $\Omega = [0,20]^2$ with the mesh sizes $\Delta x=\Delta y=20/128$. For the time-step size, we set $\Delta t=0.1$. 
The initial conditions are 
$u_1(0) = 0.001\text{rand}(x,y)$, $u_2(0) = 0.002\text{rand}(x,y)$ and $u_3(0) = 0.001\text{rand}(x,y)$.
The system parameters are taken as
$d_1 = 1$, $d_2 =2$, $d_3= 3$.

For this problem, the reduced dimension $r$ in the DMD framework is determined  by the relative cumulative energy criterion
$$ 
\min_{1\leq r \leq R}\frac{\sum_{j=1}^r \sigma_{j}^2}{\sum_{j=1}^{R} \sigma_{j}^2  } > \epsilon ,
$$
with the tolerance $\epsilon=10^{-5}$. 
The total mass is computed by the sum of averaged state vectors over the spatial domain as
$$
< \bu_i(t)  > = \frac{1}{|\Omega |} \int_\Omega \bu_i(x,y,t)dxdy.
$$

\begin{figure}[ht]
\centering
\includegraphics[width=0.32\textwidth]{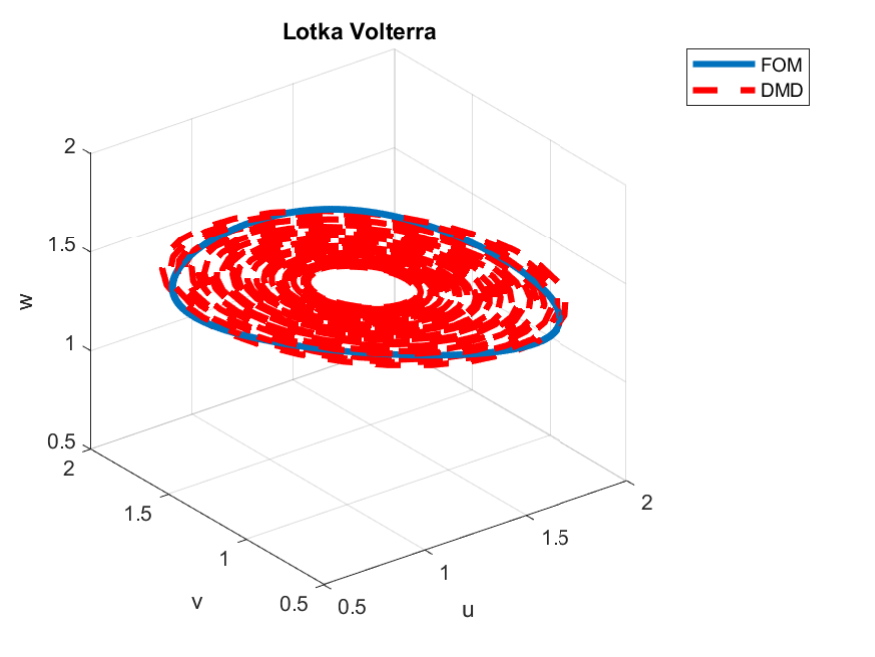}
\includegraphics[width=0.32\textwidth]{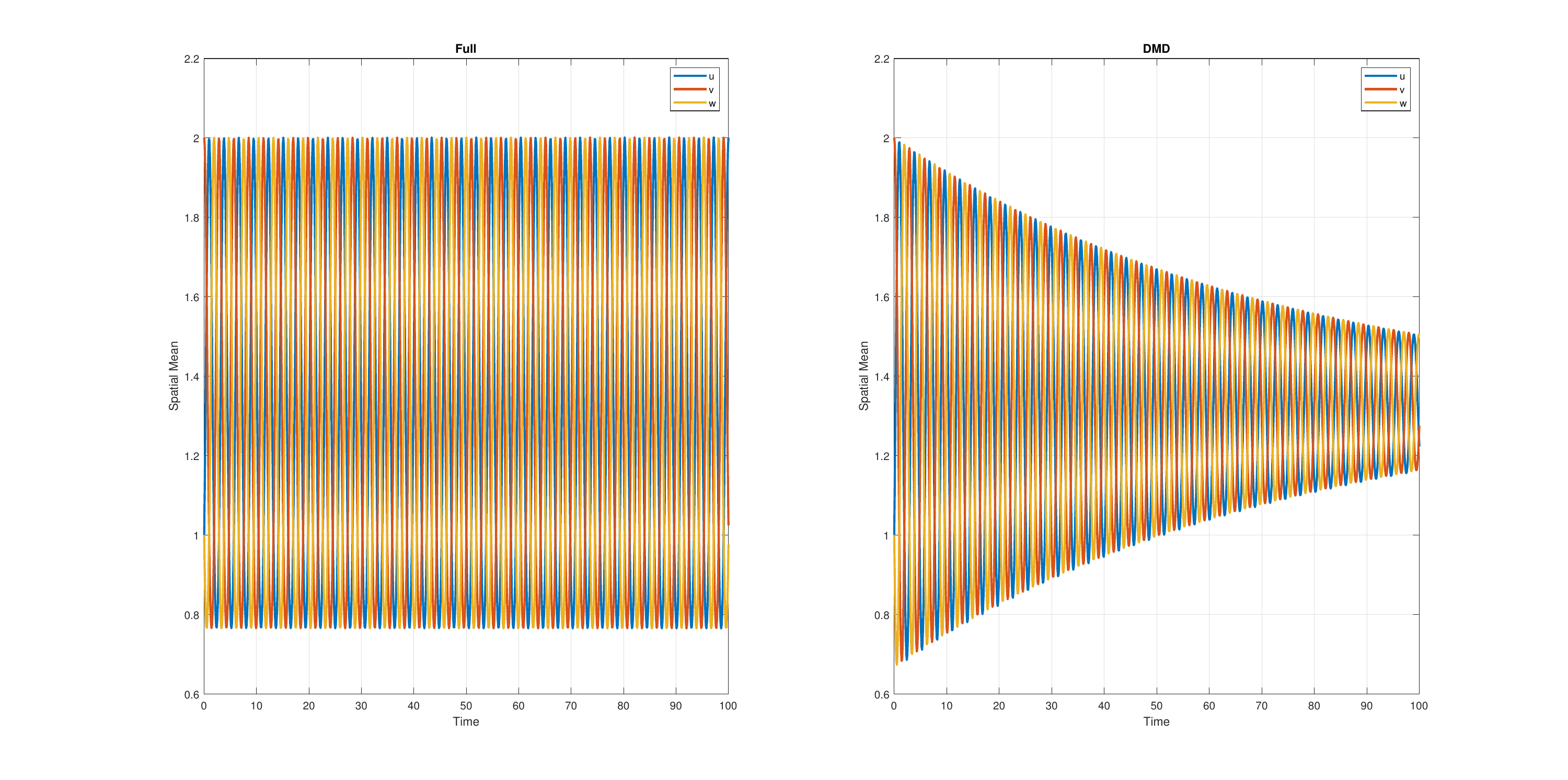}
\includegraphics[width=0.32\textwidth]{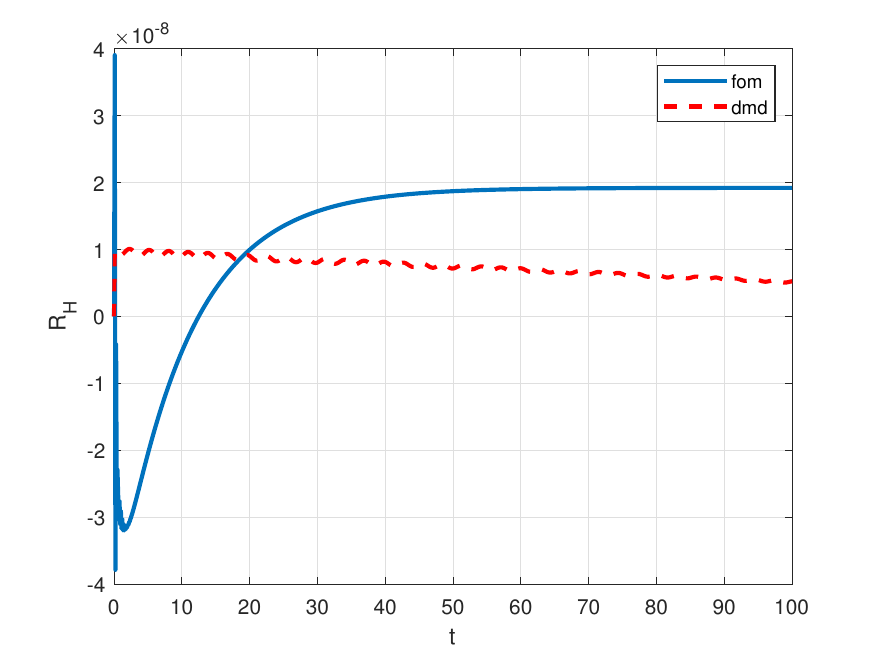}
\caption{DMD: phase portrait (left),  solutions (middle), total mass (right)}
\label{difdmd}
\end{figure}

\begin{figure}[ht]
\centering
\includegraphics[width=0.32\textwidth]{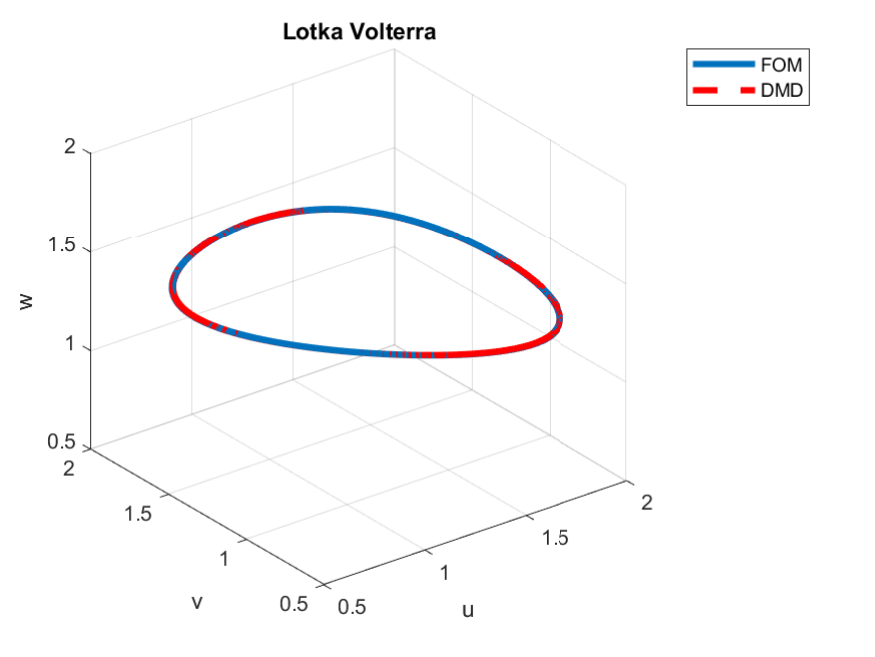}
\includegraphics[width=0.32\textwidth]{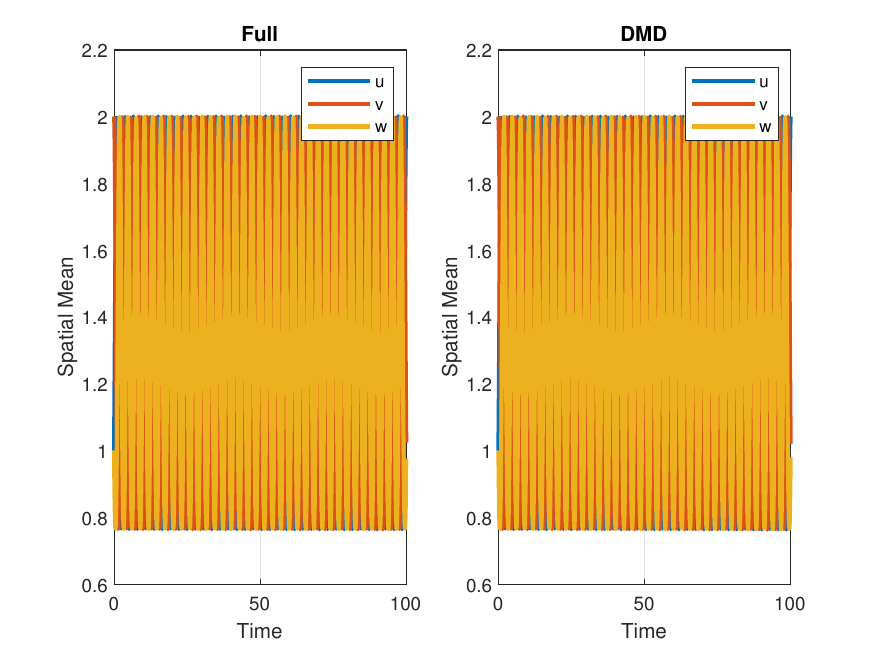}
\includegraphics[width=0.32\textwidth]{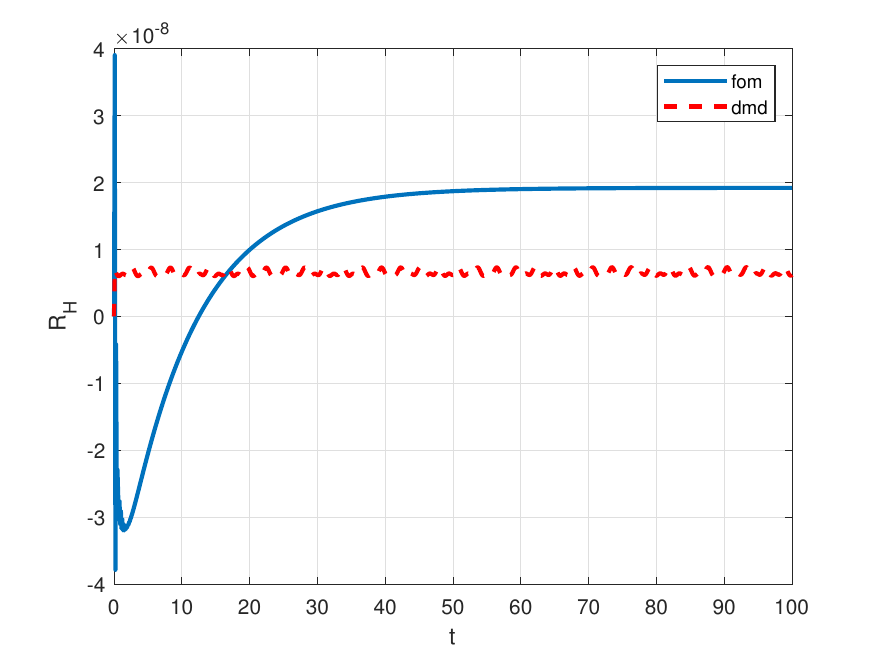}
\caption{DDMD: phase portrait (left),  solutions (middle), total mass (right)}
\label{difddmd}
\end{figure}

\begin{table}[ht]
\caption{$L_2$-errors}\label{error}
\begin{tabular*}{\textwidth}{@{\extracolsep\fill}lccccccccc}
\hline%
\multirow{2}{*}{method} & \multirow{2}{*}{rank} & \multirow{2}{*}{speedup} & \multicolumn{3}{@{}c@{}}{Training (0-50)} & & \multicolumn{3}{@{}c@{}}{Testing(50-100)} \\\cline{4-6}\cline{8-10}%
 & & & $u_1$ &  $u_2$  &  $u_3$ & &  $u_1$ &  $u_2$  &  $u_3$  \\
\hline
DMD & 2 & 13.1 &  7.15e-01 & 6.95e-01 & 7.20e-01 &  & 8.34e-01  & 8.30e-01 &  8.20e-01 \\
EDMD & 6 & 5.6 &  3.61e-06 &3.68e-06  &3.64e-06  & &  4.73e-05  & 4.60e-05  &  4.35e-05 \\
HDMD & 5 & 3.1 &  2.92e-06 & 2.74e-06 & 3.14e-06  & & 4.71e-05   & 4.58e-05 &  4.33e-05 \\
EHDMD & 7 & 0.8 & 2.63e-06  & 2.40e-06 & 2.79e-06 &  & 5.14e-05  & 4.98e-05 &  4.68e-05 \\
\hline
\end{tabular*}
\end{table}

Numerical results show that the DMD fails to predict the spatially  averaged periodic solutions, Figure \ref{difdmd} and Table \ref{error}. Meaningful results are obtained with the HDMD and EDMD as shown in Table \ref{error} and Figure \ref{difddmd}.

\section{Conclusion} 
\label{sec:conc}

HDMD and EDMD using quadratic dictionaries predict the solutions of the cyclic LVS with high accuracy. The linear Hamilton and total mass are preserved up to machine precision. The performance of the DMD solutions 
relies on the accurately simulated data by the structure preserving Kahan's integrator.


\end{document}